\documentclass{article}

\usepackage[T1]{fontenc}
\usepackage{geometry}			
\usepackage{parskip}			
\usepackage{mathtools}		
	\allowdisplaybreaks			
\usepackage{amsthm}					
\usepackage[shortlabels]{enumitem}	
\usepackage{amssymb}	

\newcommand{\card}[1]{\lvert #1 \rvert}
\newcommand{\floor}[1]{\lfloor #1 \rfloor}

\newcommand{\Set}[1]{\lbrace #1 \rbrace}
\newcommand{\defining}[1]{\emph{#1}}			

\newcommand{\union}{\cup}

\newcommand{\intersect}{\cap}
\newcommand{\Intersect}{\bigcap}
\newcommand{\bigoh}[1]{O\left(#1\right)}					
\newcommand{\bigtheta}[1]{\Theta(#1)}			
\newcommand{\littleoh}[1]{o(#1)}				

\theoremstyle{plain}
	\newtheorem*{theorem*}{Theorem}
	\newtheorem{theorem}{Theorem}[section]

	\newtheorem{proposition}[theorem]{Proposition}
	\newtheorem{conjecture}[theorem]{Conjecture}

\theoremstyle{definition}
	\newtheorem{definition}[theorem]{Definition}

	\newtheorem{example}[theorem]{Example}

\theoremstyle{remark}

\usepackage{hyperref}	
	\hypersetup{
		colorlinks=true,
		linkcolor=red,
		urlcolor=blue,
		citecolor=magenta
	}
\usepackage[abbrev,nobysame,lite]{amsrefs}

\title{Bounded fractional intersecting families are linear in size}

\author{
	Niranjan Balachandran\\
	Department of Mathematics\\
	Indian Institute of Technology Bombay\\
	Powai, Mumbai 400076\\
	Maharashtra, India\\
	Email: \texttt{niranj@math.iitb.ac.in}
	\and
	Shagnik Das\thanks{Research supported by the National Science and Technology Council (NSTC) grants 111-2115-M-002-009-MY2 and 113-2628-M-002-008-MY4.}\\
	Department of Mathematics\\
	National Taiwan University\\
	Taipei 106319\\
	Taiwan\\
	Email: \texttt{shagnik@ntu.edu.tw}
	\and
	Brahadeesh Sankarnarayanan\thanks{Corresponding author. Research supported by the National Board for Higher Mathematics (NBHM), Department of Atomic Energy (DAE), Govt.\ of India, and the Industrial Research and Consultancy Centre (IRCC), Indian Institute of Technology Bombay.}\\
	Department of Mathematics\\
	Indian Institute of Technology Bombay\\
	Powai, Mumbai 400076\\
	Maharashtra, India\\
	Email: \texttt{brahadeesh.s@iitb.ac.in}
	}

\date{May 12, 2025} 

\begin{document}
\maketitle
	
\begin{abstract}
	Using the sunflower method, we show that if \(\theta \in (0,1) \cap \mathbb{Q}\) and \(\mathcal{F}\) is a $O(n^{1/3})$-bounded \(\theta\)-intersecting family over \([n]\), then \(\lvert \mathcal{F} \rvert = O(n)\), and that if \(\mathcal{F}\) is \(o(n^{1/3})\)-bounded, then \(\lvert \mathcal{F} \rvert \leq (\frac32 + o(1))n\).
	This partially solves a conjecture of Balachandran, Mathew and Mishra that any \(\theta\)-intersecting family over \([n]\) has size at most linear in \(n\), in the regime where we have no very large sets.
\end{abstract}

\textbf{Declarations of interest}: none

\section{Introduction}\label{S:Introduction}

The study of intersecting families of set systems has a long and storied history in extremal combinatorics, with the prototypical problem taking the form: how large can a family of subsets of $[n]$ be under the constraint that the sets satisfy some intersection properties? One of the classic theorems in this direction is Fisher's Inequality~\cite{Fisher1940}, which states that if $1 \le \lambda \le n$, and $\mathcal{F} \subseteq 2^{[n]}$ is a set family with $|F \cap G| = \lambda$ for all distinct $F, G \in \mathcal{F}$, then $\card{\mathcal{F}} \le n$. With its applications to experimental design, this fundamental theorem is one of the cornerstones of design theory, where decades of work have culminated in a series of stunning existence results; by Wilson~\cite{Wilson1975} for strength-$2$ designs, and independently by Keevash~\cite{Keevash2014} and by Glock, K\"uhn, Lo, and Osthus~\cite{GKLO2023} for the general case. In the context of extremal set theory, various generalisations of Fisher's Inequality continue to provide fertile ground for research and theory-building; for instance, one can allow various different intersection sizes, place restrictions on the sizes of sets in the family, require all pairwise intersections to be the same, or consider intersections of multiple sets. Some of the highlights along these lines include the de Bruijn--Erd\H{o}s Theorem~\cite{deBruijnErdos1948}, the Erd\H{o}s--Ko--Rado Theorem~\cite{ErdosKoEtAl1961}, the Ray-Chaudhuri--Wilson Inequality~\cite{Ray-ChaudhuriWilson1975}, the Frankl--Wilson Inequality~\cite{FranklWilson1981}, the Alon--Babai--Suzuki Inequality~\cite{AlonBabaiEtAl1991}, the Erd{\H o}s--Rado Sunflower Lemma~\cite{ErdosRado1960} and the recent improvement by Alweiss, Lovett, Wu, and Zhang~\cite{AlweissLovettEtAl2021}. For a more complete survey of restricted intersection theorems, we refer the reader to the monograph of Babai and Frankl~\cite{BabaiFrankl2022}.

Returning to Fisher's Inequality, observe that it greatly restricts the set family --- we can only take a linear number of subsets from the exponentially many possibilities. However, this is perhaps not so surprising, as the condition imposed is very strong. Indeed, if $\lambda$ is large, then $\mathcal{F}$ can only consist of large sets, since we trivially require $|F| \ge \lambda$ for each $F \in \mathcal{F}$, while if $\lambda$ is small, then the large sets in our family must be essentially disjoint, and so we do not have space to pack many of them.

It is therefore natural to wonder what happens if we allow the intersections to scale with the sizes of the subsets, which leads us to the notion of fractional intersecting families, introduced by Balachandran, Mathew, and Mishra~\cite{BalachandranMathewEtAl2019}. Given \(\theta \in (0,1) \in \mathbb{Q}\), a \defining{(fractional) \(\theta\)-intersecting family} \(\mathcal{F}\) over \([n]\) is a collection of nonempty subsets of \([n]\) such that for all \(A, B \in \mathcal{F}\) with \(A \neq B\), \(\card{A \intersect B} \in \Set{ \theta\card{A}, \theta \card{B}} \).%
\footnote{More generally, given a set \(L\) of proper fractions, a \defining{(fractional) \(L\)-intersecting family} \(\mathcal{F}\) over \([n]\) is a collection of subsets of \([n]\) such that for all \(A, B \in \mathcal{F}\) with \(A \neq B\), \(\card{A \intersect B} \in \Set{\theta\card{A}, \theta \card{B}}\) for some \(\theta \in L\).}
In~\cite{BalachandranMathewEtAl2019}, the following upper bound is proved for the size of any \(\theta\)-intersecting family over \([n]\).
\begin{theorem}[Balachandran--Mathew--Mishra~\cite{BalachandranMathewEtAl2019}, 2019]
	Let $\theta = \frac{a}{b} \in (0,1) \cap \mathbb{Q}$, and let \(\mathcal{F}\) be a \(\theta\)-intersecting family over \([n]\). Then, \(\card{\mathcal{F}} = \bigoh{\frac{n \log^2 n}{\log \log n}}\), where the implicit constant depends on $b$.
\end{theorem}
On the other hand, the best-known constructions give \(\theta\)-intersecting families over \([n]\) of size only linear in \(n\).
\begin{example}\label{Eg:sunflower}
	The \defining{sunflower family} \(\mathcal{F}_{s}\) over \([n]\) is defined as follows:
	\[
		\mathcal{F}_{s} =
			\begin{cases}
				\Set{ 12, 13, \dotsc, 1n, 1234, 1256, \dotsc, 12(n-1)n}, & n \equiv 0 \pmod 2;\\
				\Set{ 12, 13, \dotsc, 1n, 1234, 1256, \dotsc, 12(n-2)(n-1)}, & n \equiv 1 \pmod 2.
			\end{cases}
	\]
	This is easily seen to be a \(\frac12\)-intersecting family, also called a \defining{bisection closed} family.
	Note that \(\card{\mathcal{F}_{s}} = \floor{3n/2} - 2\).
\end{example}

\begin{example}\label{Eg:Hadamard}
	The \defining{Hadamard family} \(\mathcal{F}_{H}\) over \([2m]\) is constructed from an \(m \times m\) normalized Hadamard matrix \(H\) as follows.
	View the rows \(A_{1},\dotsc,A_{3m}\) of the following block matrix as the \(\Set{\pm1}\)-incidence vectors of subsets of \([2m]\), where \(J\) denotes the \(m \times m\) all-ones matrix:
	\[
		\begin{bmatrix}
			H & H\\
			H & -H\\
			H & -J
		\end{bmatrix}.
	\]
	Then, \(\mathcal{F}_{H} = \Set{A_i : i \in [3m] \setminus \Set{1,m+1}}\).
	One can show using the orthogonality of the rows of \(H\) that \(\mathcal{F}_{H}\) is a bisection closed family over \([2m]\).
	Writing \(2m = n\), we see that \(\card{\mathcal{F}_{H}} = 3n/2 - 2\).
\end{example}

It was conjectured in~\cite{BalachandranMathewEtAl2019} that \emph{any} \(\theta\)-intersecting family over \([n]\) is at most linear in size.%
\footnote{The conjecture is implicit in~\cite{BalachandranMathewEtAl2019}, and explicitly stated for the case when \(\theta = 1/2\).}
\begin{conjecture}[Balachandran--Mathew--Mishra~\cite{BalachandranMathewEtAl2019}, 2019]\label{Conj:Main}
	For \(\theta \in (0,1) \intersect \mathbb{Q}\), there is a constant \(c > 0\) such that for any \(\theta\)-intersecting family \(\mathcal{F}\) over \([n]\), \(\card{\mathcal{F}} \leq cn\).
\end{conjecture}

Moreover, the fact that two very different constructions give rise to maximal bisection closed families over \([n]\) of the same size raises the question whether, for \(\theta = 1/2\), we have \(\card{\mathcal{F}} \leq \floor{3n/2} - 2\) for any bisection closed family \(\mathcal{F}\) over \([n]\).
In~\cite{BalachandranSankarnarayanan2024}, there are constructions of bisection closed families over \([n]\) for \(n \leq 15\) which have size greater than \(\floor{3n/2} - 2\), so the constructions in Examples~\ref{Eg:sunflower} and \ref{Eg:Hadamard} are possibly extremal only for large \(n\).

In this note, we make some progress towards resolving the conjecture by proving the following result.
We say that a family of sets is \defining{\(w\)-bounded}, for a positive real \(w\), if every set in the family has size at most \(w\).

\begin{theorem}\label{T:Main1}
	Let \(\theta \in (0,1) \intersect \mathbb{Q}\) and \(w = \bigoh{n^{1/3}}\) be a positive real. There is a constant \(C > 0\) such that the following holds: for all sufficiently large \(n\), if \(\mathcal{F}\) is a \(w\)-bounded \(\theta\)-intersecting family over \([n]\), then \(\card{\mathcal{F}} \leq C n\).
\end{theorem}

In fact, with the slightly stronger assumption that the family is $\littleoh{n^{1/3}}$-bounded, we can give an explicit constant that is often tight, and characterise families attaining the bound. To describe the asymptotically optimal families, we need to introduce a bit of notation related to our constructions.

First, recall that a \emph{sunflower} is a set family where all pairwise intersections are the same; that is, for all distinct $F, F' \in \mathcal{F}$, we have $F \cap F' = \bigcap_{F'' \in \mathcal{F}} F''$. The common intersection $C = \bigcap_{F'' \in \mathcal{F}} F''$ is called the \emph{core}, while the (pairwise disjoint) remainders of the sets $F \setminus C$ are called \emph{petals}. Observe that the family $\mathcal{F}_s$ from Example~\ref{Eg:sunflower} is the union of $2$- and $4$-uniform sunflowers, whose cores are nested; the definition below generalises this notion of neatly-arranged sunflowers. In what follows, given a set family $\mathcal{F}$ and a uniformity $k$, we denote by $\mathcal{F}(k)$ the collection of all sets in $\mathcal{F}$ of size $k$.

\begin{definition} \label{def:bouquet}
	Let \(\mathcal{F}\) be a family over \([n]\), and let \(k_{1} < \dotsb < k_{t}\) be the sizes of sets in \(\mathcal{F}\).
	We say that \(\mathcal{F}\) is a \emph{bouquet} if
	\begin{enumerate}
		\item each \(\mathcal{F}(k_{j})\) is a sunflower with at least two petals;
		\item\label{D2} \(C_{k_{1}} \subsetneq C_{k_{2}} \subsetneq \dotsb \subsetneq C_{k_{t}}\), where \(C_{k_{j}}\) denotes the core of \(\mathcal{F}(k_{j})\);
		\item\label{D3} for any \(F \in \mathcal{F}\) we have \(F \intersect C_{k_{t}} = C_{\card{F}}\).
	\end{enumerate}
\end{definition}

In our main result, we provide an explicit upper bound, and show that any family that attains it must essentially be a bouquet with specific set sizes.

\begin{theorem}\label{T:Main2}
	Let \(a,b \in \mathbb{N}\) such that \(1 \leq a < b\) and \(\gcd(a,b) = 1\).
	Let \(\theta = a/b\), and let \(\mathcal{F}\) be a \(\littleoh{n^{1/3}}\)-bounded \(\theta\)-intersecting family over \([n]\).
	Then \(\card{\mathcal{F}} \leq (C_{\theta} + \littleoh{1}) n\), where \(C_{\theta} = \frac{1}{b-a} \sum_{i = 1}^{\floor{b/a}} \frac{1}{i}\), and this constant is best possible for $\theta \in \Set{1/3} \union [1/2,1)$.

Furthermore, if we have equality, then there is some bouquet $\mathcal{F}^{*} \subseteq \mathcal{F}$ such that $\card{\mathcal{F} \setminus \mathcal{F}^{*}} = o(n)$, and almost all elements of $[n]$ are contained in sets of size $ib$ for all $1 \le i \le \floor{b/a}$.
\end{theorem}

\section{Main results}\label{S:Main}

Our proofs consist of two parts. First, we show that one can remove a small number of sets from a bounded \(\theta\)-intersecting family to obtain a bouquet.

\begin{proposition}\label{Prop:1}
	Let \(\theta \in (0,1) \cap \mathbb{Q}\) and \(w > 1\).
	Let \(\mathcal{F}\) be a \(w\)-bounded \(\theta\)-intersecting family over \([n]\).
	Then \(\mathcal{F}\) contains a bouquet \(\mathcal{F}^{*}\) with \(\card{\mathcal{F} \setminus \mathcal{F}^{*}} \leq w^{3}\).
\end{proposition}

We then utilise the structure of bouquets to bound their size.

\begin{proposition}\label{Prop:2}
	Let \(a,b \in \mathbb{N}\) such that \(1 \leq a < b\) and \(\gcd(a,b) = 1\).
	Let \(\theta = a/b\), and let \(\mathcal{F}^{*}\) be a \(\theta\)-intersecting bouquet over \([n]\).
	Then \(\card{\mathcal{F}^{*}} \leq C_{\theta} n\), where \(C_{\theta} = \frac{1}{b-a} \sum_{i = 1}^{\floor{b/a}} \frac{1}{i}\).
\end{proposition}

After proving these propositions, we shall combine them to deduce Theorems~\ref{T:Main1} and~\ref{T:Main2}. 

\subsubsection*{Proof of Proposition~\ref{Prop:1}}

We will require the following result of Deza~\cite{Deza1974} that implies that a large uniform \(\theta\)-intersecting family must be a sunflower.
\begin{theorem}[Deza~\cite{Deza1974}, 1974]\label{T:Deza}
	Let \(\mathcal{F}\) be a \(w\)-bounded family of subsets of \([n]\) such that all pairwise intersections have the same cardinality.
	If \(\card{\mathcal{F}} \geq w^{2} - w + 2\), then \(\mathcal{F}\) is a sunflower.
\end{theorem}

Since \(\mathcal{F}\) is \(w\)-bounded, if a level \(\mathcal{F}(k)\) is non-empty, then \(k \leq w\).
Call a level \(\mathcal{F}(k)\) \emph{small} if \(\card{\mathcal{F}(k)} \leq w^{2}\).
Note that \(\card{\mathcal{F}(1)} \leq 1 < w^{2}\): the $\theta$-intersecting property requires that any two distinct singleton sets in \(\mathcal{F}\) have intersection of size equal to $\theta \in (0,1)$, which is impossible. Thus, we can bound the number of sets in small levels by
\begin{equation}
	\sum_{k \,:\, \card{\mathcal{F}(k)} \leq w^{2}} \card{\mathcal{F}(k)} = \card{\mathcal{F}(1)} + \sum_{k > 1 \,:\, \card{\mathcal{F}(k)} \leq w^{2}} \card{\mathcal{F}(k)} \leq 1 + (w-1)w^{2} < w^{3}.
\end{equation}
We remove these sets from \(\mathcal{F}\), and shall show that what remains must be a bouquet (after removing at most one more set, if needed).
Let \(k_{1} < k_{2} < \dotsb < k_{t}\) be the remaining nonempty levels.
\begin{enumerate}
	\item By Theorem~\ref{T:Deza}, each \(\mathcal{F}(k_{j})\) is a sunflower with at least two sets.

	\item Let \(C_{k_{j}}\) be the core of \(\mathcal{F}(k_{j})\).
	Since \(\mathcal{F}\) is \(\theta\)-intersecting, \(\card{C_{k_{j}}} = \theta k_{j}\) for all \(1 \leq j \leq t\).
	Now, let \(1 \leq j < j' \leq t\), and suppose \(F' \in \mathcal{F}(k_{j'})\).
	Then \(\card{F' \intersect F} \geq \theta k_{j} = \card{C_{k_{j}}}\) for every \(F \in \mathcal{F}(k_{j})\), since \(\mathcal{F}\) is \(\theta\)-intersecting.
	If \(C_{k_{j}} \nsubseteq F'\), then \(F'\) must intersect every petal in \(\mathcal{F}(k_{j})\).
	But then \(\card{F'} \geq \card{\mathcal{F}(k_{j})} > w^{2}\), which is not possible since \(\mathcal{F}\) is \(w\)-bounded.
	Thus, \(C_{k_{j}} \subseteq F'\) for every \(F' \in \mathcal{F}(k_{j'})\), which implies that \(C_{k_{j}} \subseteq C_{k_{j'}}\).
	
	\item Let \(j < t\) and \(F \in \mathcal{F}(k_{j})\).
	If \(F \intersect (C_{k_{t}} \setminus C_{k_{j}}) \neq \emptyset\), then for any \(G \in \mathcal{F}(k_{t})\) we have \(\card{F \intersect G} > \card{C_{k_{j}}} = \theta k_{j}\).
	Thus, necessarily, \(\card{F \intersect G} = \theta k_{t}\).
	Again, \(F\) is not large enough to meet every petal of \(\mathcal{F}(k_{t})\), and so we have \(C_{k_{t}} \subseteq F\).
	Now, for any \(j' < t\), if there was another such set \(F' \in \mathcal{F}(k_{j'})\), then we would have \(\card{F \intersect F'} \geq \card{C_{k_{t}}}  = \theta k_{t} \notin \Set{\theta \card{F}, \theta \card{F'}}\), contradicting that \(\mathcal{F}\) is \(\theta\)-intersecting.
	Hence, there is at most one such set; if so, we remove it, and the remaining family satisfies \(F \cap C_{k_t} = C_{\card{F}}\).
\end{enumerate}
Then, having removed at most \(w^{3}\) sets, we are left with a bouquet \(\mathcal{F}^{*}\).\hfill\(\square\)

\subsubsection*{Proof of Proposition~\ref{Prop:2}}

Let \(k_{1} < \dotsb < k_{t}\) be the nonempty levels in the bouquet \(\mathcal{F}^{*}\) over \([n]\), and set \(Y = [n] \setminus C_{k_{t}}\).
Note that for each \(F \in \mathcal{F}^{*}\) we have \(\card{F \intersect Y} = (1 - \theta)\card{F}\).
Moreover, for each \(j\), the sets in \(\mathcal{F}^{*}(k_{j})\) are pairwise disjoint over \(Y\).
Thus, we have
\[
	\card{\mathcal{F}^{*}} = \sum_{F \in \mathcal{F}^{*}} 1 = \sum_{F \in \mathcal{F}^{*}} \sum_{y \in F \intersect Y} \frac{1}{(1-\theta)\card{F}} = \sum_{y \in Y} \sum_{\substack{F \in \mathcal{F}^{*}: \\ y \in F}} \frac{1}{(1-\theta)\card{F}}.
\]
For each \(y \in Y\), let \(S_{y} = \Set{\card{F} : F \in \mathcal{F}^{*}, y \in F}\).
Then we have
\begin{equation} \label{eqn:sizeofF*}
	\card{\mathcal{F}^{*}} = \frac{1}{1-\theta} \sum_{y \in Y} \sum_{s \in S_{y}} \frac{1}{s}.
\end{equation}
Now observe that if \(F, F' \in \mathcal{F}^{*}\) and \(\card{F} < \theta \card{F'}\), then we must have \(\card{F \intersect F'} = \theta \card{F}\).
However, \(F \intersect F' \intersect C_{i_{t}} = C_{\card{F}}\), which is of size \(\theta \card{F}\), and so \(F \intersect F' \intersect Y = \emptyset\). This means that for every \(y \in Y\), we have \(\max S_{y} \leq \frac{1}{\theta} \min S_{y}\).
Moreover, since \(\mathcal{F}^{*}\) is \(\theta\)-intersecting, \(b\) must divide \(\card{F}\) for every \(F \in \mathcal{F}^{*}\).
Thus, for every \(y \in Y\), we have some \(m_{y} \in \mathbb{N}\) such that \(S_{y} \subseteq \Set{bm_{y}, b(m_{y} + 1), \dotsc, b\floor{m_{y}/\theta}}\), and
\[
	\sum_{s \in S_{y}} \frac{1}{s} \leq \sum_{i = m_{y}}^{\floor{m_{y}/\theta}} \frac{1}{bi} = \frac{1}{b} \sum_{i = m_{y}}^{\floor{m_{y}/\theta}} \frac{1}{i}.
\]
Hence we have
\[
	\card{\mathcal{F}^{*}} = \frac{1}{1-\theta} \sum_{y \in Y} \sum_{s \in S_{y}} \frac{1}{s} \leq \frac{1}{(1-\theta)b} \sum_{y} \sum_{i = m_{y}}^{\floor{m_{y}/\theta}} \frac{1}{i}.
\]
Now write \(b = a\ell + r\), where \(\ell \in \mathbb{N}\) and \(0 \leq r \leq a - 1\).
Then,
\[
	\floor{m_{y}/\theta} = \floor{bm_{y}/a} = \ell m_{y} + \floor{rm_{y}/a} \leq \ell m_{y} + m_{y} - 1 = (\ell+1)m_{y} - 1.
\]
Thus,
\begin{align} \label{eqn:myequals1}
	\sum_{i = m_{y}}^{\floor{m_{y}/\theta}} \frac{1}{i} &= \sum_{j = 1}^{\ell - 1} \sum_{i = jm_{y}}^{(j+1)m_{y} - 1} \frac{1}{i} + \sum_{i = \ell m_{y}}^{\ell m_{y} + \floor{rm_{y}/a}} \frac{1}{i} \notag \\
	&\leq \sum_{j = 1}^{\ell -1}\sum_{i = jm_{y}}^{(j+1)m_{y} - 1} \frac{1}{jm_{y}} + \sum_{i = \ell m_{y}}^{(\ell +1)m_{y} - 1} \frac{1}{\ell m_{y}} = \sum_{j = 1}^{\ell } \frac{1}{j}. 
\end{align}
Noting that $\ell = \floor{b/a}$, \(\frac{1}{(1-\theta)b} = \frac{1}{b-a}\), and that there are at most \(n\) choices for \(y \in Y\), we obtain the desired bound.\hfill\(\square\)

Theorem~\ref{T:Main1} immediately follows from Propositions~\ref{Prop:1} and~\ref{Prop:2}. To establish Theorem~\ref{T:Main2}, we need to characterise the bouquets that attain equality in Proposition~\ref{Prop:2}.

\subsubsection*{Proof of Theorem~\ref{T:Main2}}
Let $\mathcal{F}$ be a $w$-bounded $\theta$-intersecting family over $[n]$, where $\theta = \frac{a}{b}$. By Proposition~\ref{Prop:1}, we can discard at most $w^3 = o(n)$ sets from $\mathcal{F}$ to obtain a bouquet $\mathcal{F}^* \subseteq \mathcal{F}$. By Proposition~\ref{Prop:2}, it follows that $\card{\mathcal{F}^{*}} \le C_\theta n$, where $C_{\theta} = \frac{1}{b-a} \sum_{i=1}^{\floor{b/a}} \frac{1}{i}$. Thus, we have $\card{\mathcal{F}} \le \left( C_{\theta} + o(1) \right) n$.

In order to have equality, we must have $\card{\mathcal{F}^*} \ge \left( C_\theta - o(1) \right) n$. In the proof of Proposition~\ref{Prop:2}, specifically~\eqref{eqn:sizeofF*}, we have $\card{\mathcal{F}^*} = \frac{1}{1 - \theta} \sum_{y \in Y} \sum_{s \in S_y} \frac{1}{s}$, where $Y$ is the set of elements outside the largest core of the bouquet, and $S_y$ is the set of sizes of the sets that contain $y$. In~\eqref{eqn:myequals1} we then bounded $\frac{1}{1 - \theta} \sum_{s \in S_y} \frac{1}{s}$ in terms of $m_y = \min S_y$, showing it achieves its maximum of $C_\theta$ when $S_y = \{ 1, 2, \hdots, \floor{b/a} \}$; we call this latter set $S^*$.

To prove stability, let $Y = Y_1 \cup Y_2 \cup Y_3$, where $Y_1 = \{ y \in Y: S_y = S^* \}$, $Y_2 = \{ y \in Y: m_y \ge 2\}$, and $Y_3 = \{y \in Y: m_y = 1, S_y \subsetneq S^* \}$.

For $y \in Y_3$, since $S_y$ is missing at least one element from $S^*$, we have $\frac{1}{1-\theta} \sum_{s \in S_y} \frac{1}{s} \le C_\theta - \frac{1}{(1 - \theta) \floor{b/a}}$.

For $y \in Y_2$, consider the $j = 1$ term in the inequality in~\eqref{eqn:myequals1}. Here, in the sum $\sum_{i = m_y}^{2m_y - 1} \frac{1}{i}$, we bound the summands from above by $\frac{1}{m_y}$ to show that this sum is at most $1$. However, if $m_y \ge 2$, then for the terms with $\frac32 m_y \le i \le 2m_y - 1$, of which there are at least $\frac13 m_y$, the summand is in fact at most $\frac{2}{3m_y}$. Thus, in this case, this sum is at most $\frac{8}{9}$. This shows that for $y \in Y_2$, we have $\frac{1}{1-\theta} \sum_{s \in S_y} \frac{1}{s} \le C_\theta - \frac{1}{9(1-\theta)}$.

Hence, we have
\[ \card{\mathcal{F}^*} = \frac{1}{1-\theta} \sum_{y \in Y} \sum_{s \in S_y} \frac{1}{s} \le C_\theta n - \frac{1}{9(1 - \theta)} \card{Y_2} - \frac{1}{(1-\theta)\floor{b/a}} \card{Y_3}. \]
Thus, to have $\card{\mathcal{F}^*} \ge \left( C_\theta - o(1) \right) n$, we must have $\card{Y_2} = \card{Y_3} = o(n)$, which means almost all elements of the ground set are in $Y_1$, and are thus contained in sets of size $b, 2b, \dotsc, \floor{b/a} b$.

Finally, we show this constant $C_\theta$ is best possible when $\theta \in \{ 1/3 \} \cup [1/2, 1)$ by means of the following constructions of $\theta$-intersecting $\littleoh{n^{1/3}}$-bounded families over $[n]$.
\begin{itemize}
	\item For \(\theta = a/b \in (1/2,1)\), let \(\mathcal{F}\) be a maximal \(b\)-uniform sunflower over \([n]\) with core of size \(a\).
Then \(\mathcal{F}\) is \(\frac{a}{b}\)-intersecting and \(\card{\mathcal{F}} = \floor{\frac{n-a}{b-a}} = \floor{C_{\theta}n - \frac{a}{b-a}}\).
	\item For \(\theta = 1/2\), the family \(\mathcal{F}_{s}\) has size \(\floor{3n/2} - 2\) over \([n]\), and \(C_{1/2} = 3/2\).
	\item For \(\theta = 1/3\), assume \(n \equiv 3 \pmod{24}\) for convenience, and consider the family \(\mathcal{F} = \mathcal{F}(3) \cup \mathcal{F}(6) \cup \mathcal{F}(9)\), where:
\begin{itemize}
	\item \(\mathcal{F}(3)\) is a sunflower with core \(\Set{1}\) and petals \(\Set{ \Set{2i, 2i+1} : 2 \leq i \leq (n-1)/2} \),
	\item \( \mathcal{F}(6)\) is a sunflower with core \(\Set{1,2}\) and petals \(\Set{ \Set{24i + j, 24i + j + 6, 24i + j + 12, 24i + j + 18} : 0 \leq i \leq (n-27)/24,\, 4 \leq j \leq 9}\), and
	\item \( \mathcal{F}(9)\) is a sunflower with core \(\Set{1,2,3}\) and petals \(\Set{ \Set{6i-2, 6i-1, 6i, \dotsc, 6i+3} : 1 \leq i \leq (n-3)/6}\).
\end{itemize}
	\(\mathcal{F}\) is then \(\frac13\)-intersecting, and \(\card{F} = (n-3)\left( \frac{1}{2} + \frac{1}{4} + \frac{1}{6} \right) = \frac{11}{12}(n-3) = C_{1/3}n - \frac{33}{12}\). \hfill \(\square\)
\end{itemize}

\section{Concluding remarks}

\paragraph{Sharpening the constant.} While Theorem \ref{T:Main2} establishes the correct constant for certain values of \(\theta\), further arguments can be made to improve the constant for other fractions. We briefly illustrate this with the example of \(\theta=1/4\): any \(o(n^{1/3})\)-bounded \(\frac{1}{4}\)-intersecting family has size at most \((\frac{7}{12}+ \littleoh{1})n\), rather than the $(\frac{25}{36} + o(1))n$ bound that Theorem~\ref{T:Main2} gives.
	
	For the lower bound, we construct such a family using sets of size \(4, 8\) and \(16\). The sunflowers have nested cores of size \(1, 2\) and \(4\) respectively, and  for the petals, we divide the remaining elements into blocks of size \(36\), arranged in \(3 \times 12\) rectangles. Each row (of size \(12\)) is the petal of a \(16\)-set, and is partitioned into four petals of size \(3\) each (for the \(4\)-sets). The \(12\) columns are paired up to form the petals of the \(8\)-sets, in such a way that they intersect each small petal at most once.

For the upper bound, we first note that if the constant of $\frac{25}{36}$ from the theorem were tight, then almost all elements of the ground set would have to be contained in sets of size $4$, $8$, $12$, and $16$. However, it is not hard to show (we omit the details) that sets of size $12$ are not compatible with the other set sizes. Given this, one can then show that $\sum_{s \in S_y} \frac{1}{s}$ is maximised when $S_y = \{4, 8, 16\}$, which results in a bound of $( \frac{7}{12} + o(1) )n$ instead.

It appears a difficult task to see what the correct constant is in general, even for \(\theta=1/b\), and it would be interesting to obtain further results in this direction.

\paragraph{Small families.} The \(\littleoh{n}\) error in Theorem~\ref{T:Main2} is necessary, because of the existence of bisection closed families of size greater than \(3n/2\) for \(n \leq 15\).
	These are constructed in~\cite{BalachandranSankarnarayanan2024} using the Fano plane.
	Define the family \(\mathcal{F}_{\text{Fano}}\) over \([8]\) as follows:
\begin{align*}
	\mathcal{F}_{\text{Fano}} &= \mathcal{F}_{s} \union \Set{1357,1368,1458,1467}\\
	&= \Set{12,13,14,15,16,17,18,1234,1256,1278,1357,1368,1458,1467}.
\end{align*}
It is easy to check that \(\mathcal{F}_{\text{Fano}}\) is a bisection closed family of size \(14\) over \([8]\), and it arises from the symmetric \(2\)-\((7,4,2)\) design.
We can similarly modify \(\mathcal{F}_{s}\) using the sets \(1357\), \(1368\), \(1458\), and \(1467\) to get bisection closed families over \([n]\) of size more than \(\floor{3n/2} - 2\) for \(n \leq 15\).

\paragraph{Unbounded families.} Our constant \(C_{\theta}\) in Theorem~\ref{T:Main2} is strictly smaller than 3/2 when \(\theta \neq 1/2\), and the proof of Proposition~\ref{Prop:2} shows that even when \(\theta = 1/2\), we obtain a smaller constant unless almost all elements of \([n]\) are contained in sets of size \(2\). However, the existence of the Hadamard families \(\mathcal{F}_{H}\) of Example~\ref{Eg:Hadamard} precludes any simple extension of the argument given in this note to try and establish an upper bound of \((\tfrac32 + \littleoh{1})n\) for the size of an arbitrary bisection closed family, since these are bisection closed families of size \(3n/2 - 2\) that do not contain any sets of size \(2\) (in fact, the set sizes in \(\mathcal{F}_{H}\) are all either \(n/2\) or \(n/4\)).
	
\paragraph{Families of large sets.} The best results known so far in the ``large'' regime are given in~\cite{BalachandranMathewEtAl2019}: if all the sets in \(\mathcal{F}\) have size at least \(\frac{1}{4(1-\theta)}n - \bigtheta{\sqrt{n}}\), then \(\card{\mathcal{F}} = \bigoh{n}\).

\paragraph{A linear algebraic reformulation.} In~\cite{BalachandranSankarnarayanan2024}, the authors consider a related problem of finding bounds on the ranks of certain symmetric matrices.
	Specifically, large \(\theta\)-intersecting families induce such matrices of low rank.
	There the authors construct low rank matrices using bipartite graphs and ask whether any of them arise from \(\theta\)-intersecting families.
	Theorem \ref{T:Main1}  shows that it is not possible for \emph{bounded} bisection closed families to induce such matrices.
	This explains in a sense why the Fano construction does not seem to extend beyond small values of \(n\) to produce larger bisection closed families from \(\mathcal{F}_{s}\).

\paragraph{Hierarchically closed families.} Our results also have implications in the setting of \defining{hierarchically \(r\)-closed} \(\theta\)-intersecting families, as defined by Balachandran, Bhattacharya, Kher, Mathew, and Sankarnarayanan~\cite{BalachandranBhattacharyaEtAl2023c}.
Given \(r \geq 2\), we say \(\mathcal{F}\) is hierarchically \(r\)-closed \(\theta\)-intersecting if, for any \(2 \leq t \leq r\) and any \(t\)-subset \(\Set{A_1, \dotsc, A_t}\) of \(\mathcal{F}\), we have \(\card{\Intersect_{i=1}^t A_i} \in \Set{\theta \card{A_i} : i \in [t] }\).
From our previous examples, note that \(\mathcal{F}_s\) is hierarchically \(r\)-closed for all \(r\), while \(\mathcal{F}_H\) is not hierarchically \(r\)-closed for any \(r \geq 3\).
Thus, in this sense, the two families are at opposite ends of a spectrum, despite having the same size.%
\footnote{Note that a hierarchically \(2\)-closed family is just a \(\theta\)-intersecting family as defined in Section~\ref{S:Introduction}. So, when we say that a \(\theta\)-intersecting family \(\mathcal{F}\) is \emph{hierarchically closed}, we mean that it is hierarchically \(r\)-closed for some \(r \geq 3\).}

In~\cite[Theorem 5]{BalachandranBhattacharyaEtAl2023c}, it was shown that Conjecture~\ref{Conj:Main} holds for hierarchically closed fractional intersecting families with a constant \(c_{\theta} \leq \frac{1}{b-a} (2 \log(\theta^{-1}) + 2)\). In the special case of $\theta = 1/2$, the authors improved the constant to the tight $c_{1/2} = 3/2$, and further showed that $\mathcal{F}_s$ from Example~\ref{Eg:sunflower} is the unique extremal family, up to permutation of the ground set. To prove their results, they showed that a hierarchically closed fractional intersecting family must essentially be a bouquet, and then gave an upper bound on the size of bouquets.

In Proposition~\ref{Prop:2}, we provide sharper estimates on the size of bouquets, and thus we improve the bounds in~\cite{BalachandranBhattacharyaEtAl2023c} for hierarchically closed fractional intersecting families, showing that these have size at most $( C_{\theta} + \littleoh{1})n$ as well. Furthermore, since the constructions of bouquets in Theorem~\ref{T:Main2} are also hierarchically closed, the tightness results carry over as well.

\subsection*{Acknowledgements} We would like to thank the anonymous referees for their helpful suggestions to improve the presentation of our results.


\end{document}